\documentstyle[amscd,amssymb,verbatim]{amsart}
\newcommand{\pa}{\partial}

\newcommand{\pr}{\operatorname{pr}}

\renewcommand{\Re}{\operatorname{Re}}

\newcommand{\lan}{\langle}
\newcommand{\ran}{\rangle}

\newcommand{\HS}{{\frak H}}

\newcommand{\de}{\delta}

\numberwithin{equation}{section}

\newtheorem{thm}{Theorem}
\newtheorem{prop}{Proposition}[section]
\newtheorem{lem}[prop]{Lemma}

\newtheorem{cor}[prop]{Corollary}

\newenvironment{rems}{\vspace{3mm}
\noindent {\bf Remarks.}}{\vspace{3mm}}

\newcommand{\Pf}{\noindent {\it Proof}}
\newcommand{\id}{\operatorname{id}}

\newcommand{\ov}{\overline}

\renewcommand{\Im}{\operatorname{Im}}

\newcommand{\ra}{\rightarrow}

\newcommand{\SS}{{\cal S}}

\newcommand{\OO}{{\cal O}}

\newcommand{\dbar}{\overline{\partial}}

\renewcommand{\a}{\alpha}

\newcommand{\om}{\omega}
\newcommand{\De}{\Delta}

\newcommand{\th}{\theta}
\newcommand{\C}{{\Bbb C}}

\newcommand{\Z}{{\Bbb Z}}
\newcommand{\Q}{{\Bbb Q}}

\newcommand{\sign}{\operatorname{sign}}

\newcommand{\ed}{\qed\vspace{3mm}}

\title{Rapidly convergent series for the Weierstrass zeta-function and
the Kronecker function}
\author{A. Polishchuk}
\thanks{This work was partially supported by NSF grant}

\begin{document}

\begin{abstract} We present expressions for the Weierstrass zeta-function
and related elliptic functions by rapidly convergent series. These series
arise as triple products in the $A_{\infty}$-category of an elliptic
curve.
\end{abstract}

\maketitle

\bigskip

\section{Formulas}

In this section we derive our formulas by classical means.
In the next section we'll explain how one can guess these formulas from
the computation of certain triple products on elliptic curves.

\subsection{Weierstrass zeta-function}

Let $L$ be a lattice in $\C$, $\om_1$, $\om_2$ be generators
of $L$ such that $\Im(\ov{\om_1}\om_2)>0$. The Weierstrass zeta-function
is defined by the series
$$\zeta(x,L)=\frac{1}{x}+\sum_{\om\in L-0}(\frac{1}{x+\om}-\frac{1}{\om}+
\frac{x}{\om^2})$$
for $x\in\C\setminus L$. One has
$$\zeta(x+\om,L)-\zeta(x,L)=\eta(\om)$$
for any $\om\in L$, where $\eta(\om)$ is a constant. We denote
$\eta_i=\eta(\om_i)$ for $i=1,2$. It is well-known that
\begin{equation}\label{eta}
\eta_i=2\zeta(\frac{\om_i}{2})
\end{equation}
for $i=1,2$. Following Hecke (see \cite{H}) for any
$x=x_1\om_1+x_2\om_2\in\C-L$ we set
$$Z(x,L)=\zeta(x,L)-x_1\eta_1-x_2\eta_2$$
(here $x_1$ and $x_2$ are real).

\begin{thm}\label{mainthm} 
For any $x\in\C-L$ one has the following identity
\begin{equation}\label{mainid}
Z(x,L)=\sum_{\om\in L}
\frac{\exp(-\frac{\pi}{a(L)}|\om+x|^2)}{\om+x}-
\sum_{\om\in L-0} 
\frac{\exp(-\frac{\pi}{a(L)}|\om|^2+2\pi i E_L(\om,x))}{\om}
\end{equation}
where $a(L)=\Im(\ov{\om_1}\om_2)$ is the area of $\C/L$,
$$E_L(x,y)=\frac{\Im(\ov{x}y)}{a(L)}=\frac{\ov{x}y-x\ov{y}}{2ia(L)}$$
is the symplectic form on $\C$ (considered as a real space)
associated with the oriented lattice $L$.
\end{thm}

Let $\SS(\C)$ be the Schwarz space of $\C$. For any $\varphi\in\SS(\C)$ we
define its {\it symplectic} Fourier transform by the formula
$$\hat{\varphi}(y)=\int_{x\in\C}\varphi(x)\exp(2\pi i E_L(y,x))d_Lx$$
where $d_Lx$ is the Haar measure on $\C$ normalized by the
condition $\int_{\C/L}d_Lx=1$.

We will need the following simple lemma.

\begin{lem}\label{Poisson}
For any $\varphi\in\SS(\C)$, any $x,y\in\C$ one has
$$\sum_{\om\in L}\varphi(\om+x)\exp(-2\pi i E_L(\om+x,y))=
\sum_{\om\in L}\hat{\varphi}(\om+y)
\exp(-2\pi i E_L(\om,x))$$
\end{lem}

\Pf . By Poincare summation formula the distribution
$\de_L=\sum_{\om\in L}\de_{\om}$ is Fourier self-dual.
Since the translation by $x$ goes under Fourier transform
to the multiplication by $\exp(-2\pi i E_L(?,x))$ the result
follows.
\ed

\noindent
{\it Proof of theorem \ref{mainthm}}.
Let us denote
$$f(x)=\sum_{\om\in L}
\frac{\exp(-\frac{\pi}{a(L)}|\om+x|^2)}{\om+x}-
\sum_{\om\in L-0} 
\frac{\exp(-\frac{\pi}{a(L)}|\om|^2+2\pi i E_L(\om,x))}{\om}$$
where $x\in\C-L$.
We claim  that
$$\ov{\pa}f(x)=-\frac{\pi}{a(L)}.$$
Indeed, we have
$$\ov{\pa}f(x)=-\frac{\pi}{a(L)}\cdot\left(\sum_{\om\in L}
\exp(-\frac{\pi}{a(L)}|\om+x|^2)-\sum_{\om\in L-0}
\exp(-\frac{\pi}{a(L)}|\om|^2+2\pi i E_L(\om,x))\right),$$
so our claim follows from Lemma \ref{Poisson} since the function
$\exp(-\frac{\pi}{a(L)}|x|^2)$ goes to itself under the
symplectic Fourier transform.

It follows that the function $g(x)=f(x)+\frac{\pi}{a(L)}\ov{x}$
is holomorphic on $\C-L$.
Furthermore, looking at the series for $f$ we immediately see that
$g$ has simple poles at all the lattice points $\om\in L$ with
residues equal to $1$.

On the other hand, from the fact that the
symplectic form $E_L$ takes integer values on $L$ one immediately
derives that $f(x+\om)=f(x)$ for any $\om\in L$. Thus, we have
$$g(x+\om)=g(x)+\frac{\pi}{a(L)}\ov{\om}$$
for $\om\in L$. The Legendre period relation
$$\eta_1\om_2-\eta_2\om_1=2\pi i$$
implies that there exists a constant $c\in\C$ such
that
\begin{equation}\label{const}
\eta_i=c\om_i+\frac{\pi}{a(L)}\ov{\om_i}
\end{equation}
for $i=1,2$.
It follows that $h(x)=g(x)-\zeta(x,L)+cx$ 
is a holomorphic function on $C-L$, periodic with respect to $L$.
Comparing the polar parts of $g$ and $\zeta$ at the lattice points
we conclude that $h$ is holomorphic on $\C$. Therefore, $h$ is
constant and we have
$$f(x)-\zeta(x,L)=-\frac{\pi}{a(L)}\ov{x}-cx+h.$$
From the definition of the constant $c$ we derive that
$$x_1\eta_1+x_2\eta_2=cx+\frac{\pi}{a(L)}\ov{x}$$
where $x=x_1\om_1+x_2\om_2$.
Thus, we have
$$f(x)-\zeta(x,L)=-x_1\eta_1-x_2\eta_2+h.$$
It is easy to see that if $x\in\frac{1}{2}L-L$ then
$f(x)$. Thus, substituting $x=\frac{\om_1}{2}$
and using the identity (\ref{eta}) we derive that $h=0$.
\ed

\begin{rems} 1.
It is not obvious that the right hand side of (\ref{mainid})
is holomorphic in $\om_1$, $\om_2$. This fact is equivalent to the
identity
$$\sum_{\om\in L}
(\om+x)\exp(-\frac{\pi}{a(L)}|\om+x|^2)=
\sum_{\om\in L} 
\om\exp(-\frac{\pi}{a(L)}|\om|^2+2\pi i E_L(\om,x))$$
which can be easily deduced from Lemma \ref{Poisson}.

\noindent
2. It is well-known that the following Epstein's zeta function
$$\varphi_1(s,L,x)=\sum_{\om\in L}\frac{1}{(\om+x)|\om+x|^{2s-1}}$$
defined for $\Re(s)>1$
extends to an entire function of $s$ (for fixed $x\in\C-L$).
It was shown by N.~Katz (Cor. 3.2.24 of \cite{K}) that
$$\varphi_1(\frac{1}{2},L,x)=Z(x,L)$$
for $x\in\Q L-L$.
It is not clear whether there exists an expression for 
$\varphi_1(s,L,x)$ for arbitrary $s$ similar to the one in (\ref{mainid}).
\end{rems}

Differentiating the identity (\ref{mainid}) we obtain
the following series for the Weierstrass $\wp$-function
$\wp(x)=-\zeta'(x)$:
\begin{equation}
\wp(x)=-c+
\sum_{\om\in L}
\frac{(1+\frac{\pi}{a(L)}|\om+x|^2)\exp(-\frac{\pi}{a(L)}|\om+x|^2)}
{(\om+x)^2}+
\frac{\pi}{a(L)}\cdot\sum_{\om\in L-0} 
\frac{|\om|^2\exp(-\frac{\pi}{a(L)}|\om|^2+2\pi i E_L(\om,x))}{\om^2}
\end{equation}
where the constant $c$ is determined from (\ref{const}).
Differentiating one more time we get
\begin{equation}
\wp'(x)=-\sum_{\om\in L}
\frac{(1+(1+\frac{\pi}{a(L)}|\om+x|^2)^2)\exp(-\frac{\pi}{a(L)}|\om+x|^2)}
{(\om+x)^3}+
\frac{\pi^2}{a(L)^2}\cdot\sum_{\om\in L-0} 
\frac{|\om|^4\exp(-\frac{\pi}{a(L)}|\om|^2+2\pi i E_L(\om,x))}{\om^3}.
\end{equation}

\subsection{Kronecker function}

Let us consider the following holomorphic function in $3$
variables $\tau$, $x$, $y$, where $\Im(\tau)>0$,
$0<\Im(x),\Im(y)<\Im(\tau)$:
$$F(x,y;\tau)=-\sum_{(m+\frac{1}{2})(n+\frac{1}{2})>0}
\sign(m+\frac{1}{2})\exp(2\pi i mn\tau+2\pi imx+2\pi iny)$$
where $m,n$ are integers (our choice of sign is compatible with
the notation in Zagier's paper \cite{Z}, but our variables $x$ and $y$
differ from those used in \cite{Z} by the factor $2\pi i$).
We call it the Kronecker function since Kronecker discovered (see \cite{Kr})
the following remarkable identity:
\begin{equation}\label{Kr-id}
F(x,y;\tau)= \frac{\th'_{11}(0,\tau)}{2\pi i}\cdot\frac{\th_{11}(x+y,\tau)}
{\th_{11}(x,\tau)\th_{11}(y,\tau)}
\end{equation}
where
$$\th_{11}(x,\tau)=\sum_{n\in\Z}(-1)^n\exp(\pi i (n+\frac{1}{2})^2\tau+
2\pi i (n+\frac{1}{2})x),$$
$\th'_{11}$ is the derivative of $\th_{11}(x,\tau)$ with respect to $x$.
In particular, this identity gives a meromorphic continuation of
$F$ to $\HS\times\C^2$ with poles along the divisors
$x\in L_{\tau}$, $y\in L_{\tau}$, where $L_{\tau}=\Z+\Z\tau$.

\begin{thm} One has the following identity
\begin{equation}\label{Kr-ser}
\begin{array}{l}
2\pi i F(x,y,\tau)=
\exp(-\frac{\pi}{\Im\tau}x(y-\ov{y}))\cdot \sum_{\om\in L_{\tau}}
\frac{\exp(-\frac{\pi}{\Im\tau}|\om+x|^2-2\pi i E(\om,y))}{\om+x}+\\
\exp(-\frac{\pi}{\Im\tau}y(x-\ov{x}))\cdot \sum_{\om\in L_{\tau}}
\frac{\exp(-\frac{\pi}{\Im\tau}|\om+y|^2-2\pi i E(\om,x))}{\om+y}
\end{array}
\end{equation}
where $E=E_{L_{\tau}}$.
\end{thm}

\Pf . For a fixed $\tau$ let $f(x,y)$ be the function in the
right hand side of (\ref{Kr-ser}).
First we have to check that $f(x,y)$ is meromorphic in $x$ and $y$.
Let us denote $a=\Im(\tau)$. We have
\begin{align*}
&\frac{\pa}{\pa \ov{x}}f=
-\frac{\pi}{a}\cdot
\exp(-\frac{\pi}{a}x(y-\ov{y}))\cdot \sum_{\om\in L_{\tau}}
\exp(-\frac{\pi}{a}|\om+x|^2-2\pi i E_L(\om,y))+\\
&\frac{\pi}{a}\cdot\exp(-\frac{\pi}{a}y(x-\ov{x}))\cdot
\sum_{\om\in L_{\tau}}
\exp(-\frac{\pi}{a}|\om+y|^2-2\pi i E_L(\om,x)).
\end{align*}
Thus, we have to prove that
$$\sum_{\om\in L_{\tau}}
\exp(-\frac{\pi}{a}|\om+x|^2-2\pi i E_L(\om+x,y))=
\sum_{\om\in L_{\tau}}
\exp(-\frac{\pi}{a}|\om+y|^2-2\pi i E_L(\om,x)).$$
But this follows easily from Lemma \ref{Poisson} since the
function $\exp(-\frac{\pi}{a}|x|^2)$ is Fourier self-dual.

Next we observe that for any $\om\in L$ one has
$$f(x+\om,y)=\exp(-\frac{\pi}{a}(\om-\ov{\om})y)f(x,y).$$
Since $f(x,y)=f(y,x)$ we conclude that $f$ has the same
quasi-periodicity equations as $F$. Hence,
$f/F$ is periodic with respect to $L_{\tau}$ in both variables.
Since both $f$ and $F$ have poles of the first order at
$x\in L_{\tau}$ and $y\in L_{\tau}$ the only possible poles of
$f/F$ can come from zeroes of $F$. But the only zeroes of $F$
are the zeroes of the first order along the divisor $x+y\in L_{\tau}$.
On the other hand, one immediately checks that $f(x,-x)=0$.
Therefore, $f/F$ is holomorphic, so it should be constant.
Now the identity follows by the comparison of the residues of $F$ and
$f$ at $x=0$.
\ed

From identities (\ref{mainid}) and (\ref{Kr-ser}) one immediately
deduces the following result.

\begin{cor} One has
$$\left(2\pi i F(x,y,\tau)-\frac{1}{y}\right)|_{y=0}=\zeta(x,L_{\tau})-
x\eta_1$$
where as generators of $L_{\tau}$ we take $\om_1=1$, $\om_2=\tau$.
\end{cor}

Another way to express the relation between the function $F$ and
Weierstrass zeta-function is the following:
$$\pi i\cdot (F(x,y,\tau)+F(x,-y,\tau))|_{y=0}=
\zeta(x,L_{\tau})-x\eta_1.$$
This can be deduced either from the above corollary or using
(\ref{Kr-id}) and the formula
$$\frac{\th'_{11}(x)}{\th_{11}(x)}=\zeta(x,L_{\tau})-x\eta_1$$
which can be seen from the decomposition of $\th_{11}$
into an infinite product.

\section{Explanation}

Both the functions $Z(x,L)$ and $F(x,y,\tau)$ have nice modular
properties. The modular forms $Z(x_1\om_1+x_2\om_2,L)$
for fixed $x_1,x_2\in\Q$ were considered by Hecke in \cite{H} 
(he called them ``Teilwerte'' of the Weierstrass zeta-function). 
The modular equation for $F(x,y,\tau)$ can be found in \cite{Z}.
The formulas (\ref{mainid}) and (\ref{Kr-ser}) provide an
alternative explanation of modularity but each of the series in the right
hand side is non-holomorphic (only the difference of two such series is). 
In this section we
show that this is related to the computation of certain triple products
on elliptic curve using non-holomorphic data (namely, hermitian metrics).
However, since the result doesn't depend on a choice of non-holomorphic data
the resulting expressions are holomorphic in the modular parameter.

We refer to \cite{P-ainf}
for the general discussion of higher products on elliptic curve.
The triple products related to the
two series considered in the previous section are of the following
type. Let $L$, $M_1$, $M_2$ be hermitian line bundles of degree $1$ on
a complex elliptic curve. Then one can consider a triple product
$$m_3:H^0(E,M_1)\otimes H^1(E,L^{-1})\otimes H^0(E,M_2)\ra
H^0(E,M_1M_2L^{-1}).$$
Recall that it is given by the formula
\begin{equation}\label{m3}
m_3(s_1,e,s_2)=\pr(Q(s_1e)s_2-s_1 Q(es_2))
\end{equation}
where $s_i\in H^0(E,M_i)$, $i=1,2$ the class $e\in H^1(E,L^{-1})$
is represented by a harmonic $(0,1)$-form, $Q=\dbar^* G_{\dbar}$
where $G_{\dbar}$ is the Green operator corresponding to the
laplacian $\De_{\dbar}=\dbar\dbar^*+\dbar^*\dbar$ (where $\dbar^*$
is conjugate to $\dbar$ with respect to the hermitian metric),
$\pr=\id-G_{\dbar}\De_{\dbar}$ is the harmonic projector.
The formula (\ref{m3}) shows that in fact our triple product depends
only on the operator $Q$ acting on forms with values in $M_1L^{-1}$
and $M_2L^{-1}$. Both these line bundles are of degree zero so 
(up to switching $M_1$ and $M_2$) the following three possibilities
can occur:

\noindent
(a) $M_iL^{-1}\not\simeq\OO_E$ for $i=1,2$. Then we have $Q=(\dbar)^{-1}$
in the formula (\ref{m3}) so this triple product doesn't depend
on metrics. We will show that in this case $m_3$ is expressed in
terms of the series appearing in the formula (\ref{Kr-ser}).

\noindent
(b) $M_1L^{-1}\simeq\OO_E$, $M_2L^{-1}\not\simeq\OO_E$. In this
case the operator $Q$ depends on a hermitian metric on $\OO_E$.
However, there is a natural choice of a constant metric on $\OO_E$
and it is easy to see that $Q$ doesn't change if we rescale
a metric by a constant. In this case we'll express $m_3$ in terms
of the series from the formula (\ref{mainid}).

\noindent
(c) $M_1L^{-1}\simeq M_2L^{-1}\simeq\OO_E$. One can easily see
that in this case $m_3=0$.

To compute triple products in the cases (a) and (b) we represent
our elliptic curve in the form $E=\C/\Z+\Z\tau$. We choose
$L$ to be the line bundle on $E$ such that the theta-function
$$\th(z)=\th(z,\tau)=\sum_{n\in\Z}\exp(\pi i\tau n^2+2\pi i n z)$$
descends to a section of $L$. Thus, the pull-back of $L$ to
$\C$ is canonically trivialized. For $u\in\C$ let us denote by $L(u)$
the line bundle $t_u^*L$ where $t_u:E\ra E$ is the translation by
$u$ (note that a choice of $u\in\C$ induces a trivialization of the
pull-back of $L(u)$ to $\C$).
We define the hermitian metric
on $L(u)$ by the formula
$$\lan f,g\ran_{L(u)}=\int_{C/\Z+\Z\tau}f(x)\ov{g(x)}
\exp(-2\pi a (x_2^2+2x_2u_2)dx_1dx_2$$
where we use real coordinates $x_1,x_2$ defined by
$x=x_1+x_2\tau$, so that $u=u_1+u_2\tau$ and we denote $a=\Im\tau$.
With respect to this metric one has
$$||t_u^*\th||^2=\frac{1}{\sqrt{2a}}\exp(2\pi a u_2^2).$$
As line bundles $M_1$ and $M_2$ we take $L(u)$, $L(v)$ for
some $u,v\in\C$, so that we have natural choice of sections
$s_1=t_u^*\th$, $s_2=t_v^*\th$. As a harmonic $(0,1)$-form representing
a non-trivial class $e\in H^1(E,L^{-1})$ we take
$$\a=\frac{\pi\sqrt{2}}{\sqrt{a}}\ov{\th(x)}\exp(-2\pi a x_2^2)d\ov{x}.$$

Our computation will be based on the following formula which was proven
in \cite{P-ainf} (it is equivalent to eq. (2.2.1) of \cite{P-ainf}, on the
other hand, it can be deduced from Prop. 4.1 of \cite{Shim}):
\begin{equation}\label{thpr}
\begin{array}{l}
\th(x+y)\ov{\th(x+z)}\exp(-2\pi a(x_2^2+2x_2z_2))=\\
\frac{1}{\sqrt{2a}}\sum_{m,n}(-1)^{mn}\exp(-\frac{\pi}{2a}(|m\tau-n|^2+
2(m\ov{\tau}-n)y-2(m\tau-n)\ov{z}+(y-\ov{z})^2))\varphi_{y-z,m,n}(x)
\end{array}
\end{equation}
where we denote
$$\varphi_{w,m,n}(x)=\exp(2\pi i(mx_1+(n-w)x_2)),$$
the summation is over $(m,n)\in\Z^2$.
Note that $(\varphi_{w,m,n})_{(m,n)\in\Z^2}$ descend to the orthonormal basis
of sections on $L(w)L^{-1}$. We have
$$\dbar\varphi_{w,m,n}=\frac{\pi}{a}(m\tau-n+w)\varphi_{w,m,n}d\ov{x}.$$
In the case $w\not\in\Z+\Z\tau$ this allows us to
compute $Q=(\dbar)^{-1}$ in terms of coefficients with the basis
$\varphi_{w,m,n}d\ov{x}$. In the case $w=0$ the operator $Q$ still coincides
with $(\dbar)^{-1}$ on $\varphi_{0,m,n}d\ov{x}$ for $(m,n)\neq (0,0)$.
On the other hand, $\varphi_{0,0,0}=1$ and we have
$$\dbar^*(d\ov{x})=0,$$
hence, $Q(d\ov{x})=0$.

Let us first consider the case (a). Then 
we have
$$m_3(t_u^*\th,\a,t_v^*\th)=\pr(h_u t_v^*\th-h_v t_u^*\th)$$
where we denote $h_w=Q(t_w^*\th\cdot\a)$. Using formula (\ref{thpr})
we get that for any $w\not\in\Z+\Z\tau$ one has
$$h_w=\sum_{m,n}a_{m,n}(w)\varphi_{w,m,n}$$
where
\begin{equation}\label{amn}
a_{m,n}(w)=
\frac{(-1)^{mn}\exp(-\frac{\pi}{2a}(|m\tau-n|^2+
2(m\ov{\tau}-n)w+w^2))}{m\tau-n+w}.
\end{equation}
By definition the above triple product is proportional to
$t_{u+v}^*\th$ so the computation reduces to calculating
the coefficient
\begin{equation}\label{coeff}
\frac{\lan m_3(t_u^*\th,\a,t_v^*\th),t_{u+v}^*\th\ran}
{||t_{u+v}^*\th||^2}
=\sqrt{2a}\exp(-2\pi a(u_2+v_2)^2)\cdot
\left(\lan h_u t_v^*\th,t_{u+v}^*\th\ran-
\lan h_v t_u^*\th,t_{u+v}^*\th\ran\right).
\end{equation}
Now by definition we have
\begin{align*}
&\lan h_u t_v^*\th,t_{u+v}^*\th\ran=
\int_{\C/\Z+\Z\tau} h_u(x)\th(x+v)\ov{\th(x+u+v)}
\exp(-2\pi a(x_2^2+2x_2(u_2+v_2)))dx_1dx_2=\\
&\lan h_u, t_{u+v}^*\th\cdot \ov{t_v^*\th}\exp(-2\pi a(x_2^2+2x_2v_2))\ran
\end{align*}
where the last scalar product is taken with respect to the metric
on $L(u)L^{-1}$. Applying (\ref{thpr}) we get
$$t_{u+v}^*\th\cdot \ov{t_v^*\th}\exp(-2\pi a(x_2^2+2x_2v_2))=
\sum_{m,n} b_{m,n}(u,v)\varphi_{u,m,n}$$
where
\begin{equation}\label{bmn}
b_{m,n}(u,v)=
\frac{1}{\sqrt{2a}}(-1)^{mn}\exp(-\frac{\pi}{2a}(|m\tau-n|^2+
2(m\ov{\tau}-n)(u+v)-2(m\tau-n)\ov{v}+(u+v-\ov{v})^2)).
\end{equation}
Since $\varphi_{u,m,n}$ is an orthonormal system we derive
$$\lan h_u t_v^*\th,t_{u+v}^*\th\ran=\sum_{m,n}a_{m,n}(u)
\ov{b_{m,n}(u,v)}.$$
Substituting the expressions (\ref{amn}) and (\ref{bmn}) and
simplifying we obtain
$$\sqrt{2a}\exp(-2\pi a(u_2+v_2)^2)\cdot
\lan h_u t_v^*\th,t_{u+v}^*\th\ran=
\exp(\frac{\pi}{a}u(v-\ov{v}))\cdot\sum_{\om\in\Z+\Z\tau}
\frac{\exp(-\frac{\pi}{a}|\om+u|^2+2\pi i E(\om,v))}{\om+u}
$$
where $E$ is the symplectic form associated with the oriented lattice
$\Z+\Z\tau$. Substituting this into the equation (\ref{coeff}) we obtain
that the coefficient of $m_3(t_u^*\th,\a,t_v^*\th)$ with
$t_{u+v}^*\th$ is equal to
$$\exp(\frac{\pi}{a}u(v-\ov{v}))\cdot\sum_{\om\in\Z+\Z\tau}
\frac{\exp(-\frac{\pi}{a}|\om+u|^2+2\pi i E(\om,v))}{\om+u}-
\exp(\frac{\pi}{a}v(u-\ov{u}))\cdot\sum_{\om\in\Z+\Z\tau}
\frac{\exp(-\frac{\pi}{a}|\om+v|^2+2\pi i E(\om,u))}{\om+v}.
$$
According to the identity (\ref{Kr-ser}) this expression is
equal to $2\pi i F(u,-v,\tau)$. In \cite{P-ainf} it was shown
that the series defining $F(u,-v,\tau)$ appears as the corresponding
triple product in the Fukaya category of the torus, so we can view
the identity (\ref{Kr-ser}) as a manifestation of the homological
mirror symmetry\footnote{In \cite{P-hmc} we worked with the RHS
of the identity (\ref{Kr-ser}) without presenting an explicit series
for it.}. 

In the case (b) we set $u=0$ and slightly modify the above computation.
Namely, we have
$$h_0=Q(\th\cdot\a)=\sum_{(m,n)\neq(0,0)}a_{m,n}(0)\varphi_{0,m,n}$$
where $a_{m,n}(0)$ are still given by the formula (\ref{amn}).
The formula (\ref{coeff}) still holds for $u=0$, so we obtain that the
coefficient of the triple product in this case is equal to
$$\sum_{\om\in\Z+\Z\tau,\om\neq 0}
\frac{\exp(-\frac{\pi}{a}|\om|^2+2\pi i E(\om,v))}{\om}-
\sum_{\om\in\Z+\Z\tau}
\frac{\exp(-\frac{\pi}{a}|\om+v|^2)}{\om+v}.
$$
According to the identity (\ref{mainid}) this is equal to $-Z(v,\Z+\Z\tau)$.
Thus, considering triple products above one can discover
the identity (\ref{mainid}) as follows. One starts with the identity 
(\ref{Kr-ser}) which follows from the homological mirror symmetry for 
elliptic curve. Then one tries to pass to the limit as $u\ra 0$.
From the expression of the Weierstrass zeta-function as the logarithmic
derivative of sigma-function it is easy to guess that the limit of 
$F(u,-v)$ should be related to $\zeta(v,\Z+\Z\tau)$. 
On the other hand, it should be modular,
so one naturally arrives to considering $Z(v,\Z+\Z\tau)$.

\end{document}